\documentclass[a4paper]{amsart}
\usepackage{amssymb}
\usepackage[all]{xy}
\usepackage{amsmath}
\usepackage{graphics}
\CompileMatrices

\newtheorem{theorem}{Theorem}[section]

\newtheorem{proposition}[theorem]{Proposition}
\newtheorem{corollary}[theorem]{Corollary}
\theoremstyle{definition}

\newtheorem{example}[theorem]{Example}

\newtheorem{remark}[theorem]{Remark}

\numberwithin{equation}{theorem}
\def\mal{\! \cdot \!}

\def\bangle#1{\langle #1 \rangle}

\def\KK{{\mathbb K}}

\def\ZZ{{\mathbb Z}}

\def\QQ{{\mathbb Q}}

\def\Spec{{\rm Spec}}
\def\Proj{{\rm Proj}}

\def\cone{{\rm cone}}

\newcounter{itemnumber}

\begin{document}

\sloppy

\title[On the multiplication map]
{On the multiplication map 
\\
of a multigraded algebra}

\author[I.~Arzhantsev]{Ivan~V. Arzhantsev} 
\thanks{Supported by INTAS YS 05-109-4958}
\address{Department of Higher Algebra, 
Faculty of Mechanics and Mathematics, 
Moscow State Lomonosov University,
Vorobievy Gory, GSP-2, Moscow, 119992, Russia}
\email{arjantse@mccme.ru}
\author[J.~Hausen]{J\"urgen Hausen} 
\address{Mathematisches Institut, Universit\"at T\"ubingen,
Auf der Morgenstelle 10, 72076 T\"ubingen, Germany}
\email{hausen@mail.mathematik.uni-tuebingen.de}
\subjclass[2000]{13A02, 14L24}

\begin{abstract}
Given a multigraded algebra $A$,
it is a natural question whether
or not for two homogeneous components
$A_u$ and $A_v$, the product
$A_{nu}A_{nv}$ is the whole component 
$A_{nu+nv}$ for $n$ big enough.
We give combinatorial and geometric
answers to this question.
\end{abstract}

\maketitle

\section{Statement and discussion of the results}

In this note, we consider the multiplication
map of a multigraded algebra 
and ask for its surjectivity properties 
on the homogeneous parts.
More precisely, 
let $A$ be an (associative, commutative), 
integral, finitely generated algebra 
(with unit) over an algebraically closed field~$\KK$,
and suppose that $A$ is graded by a lattice
$M \cong \ZZ^d$, i.e., we have
\begin{eqnarray*}
A
& = &
\bigoplus_{u \in M} A_u.
\end{eqnarray*}
By the weight cone of $A$
we mean the convex, polyhedral cone
$\omega(A) \subseteq \QQ \otimes_\ZZ M$ 
generated by all 
$u \in M$ with $A_u \ne 0$.
We investigate the following
problem:
given $u,v \in \omega(A) \cap M$,
does there exist an $m > 0$ such that for any 
$k > 0$ the multiplication map
defines a surjection
$$ 
\mu_{km} \colon 
A_{kmu} \otimes_\KK A_{kmv}
\ \to \ 
A_{km(u+v)},
\qquad
f \otimes g
\ \mapsto \
fg.
$$
We call a pair $u,v \in \omega(A) \cap M$ 
{\em generating\/} if it has this property.
Simple examples show that not every 
pair is generating.
In our first result 
we provide combinatorial criteria for 
a pair to be generating,
and in the second one, we give a
geometric characterization  
for the case of a factorial algebra
$A$.

To present the first result,
let us  recall from~\cite{BeHa2}
the concept of the GIT-fan 
associated to~$A$.
The $M$-grading of $A$ defines a
(unique) action of the torus 
$T := \Spec(\KK[M])$
on $X := \Spec(A)$ such that
for any $u \in M$, the elements 
$f \in A_u$ are precisely the 
semiinvariants of the 
character $\chi^u \colon T \to \KK^*$,
i.e., each $f \in A_u$ satisfies
\begin{eqnarray*}
f(t \mal x)
& := & 
\chi^u(t) f(x).
\end{eqnarray*}
The {\em orbit cone\/} of a (closed) point 
$x \in X$ is the convex, polyhedral cone 
$\omega(x) \subseteq \QQ \otimes_\ZZ M$
generated by all $u \in \omega(A)$ admitting
an $f \in A_u$ with $f(x) \ne 0$.
The collection of orbit cones 
is finite, and thus one may associate to
any element $u \in \omega(A)$ its,
again convex, polyhedral, {\em GIT-cone}:
\begin{eqnarray*}
\lambda(u)
& := & 
\bigcap_{\substack{x \in X, \\ u \in \omega(x)}} \omega(x).
\end{eqnarray*}
These GIT-cones cover the weight cone 
$\omega(A)$, and by~\cite[Thm.~3.11]{BeHa2}, 
the collection $\Lambda(A)$ of all of them
is a fan in the sense that
if $\lambda \in \Lambda(A)$ then 
also every face of $\lambda$ belongs to $\Lambda(A)$,
and for $\tau,\lambda \in \Lambda(A)$, 
the intersection
$\tau \cap \lambda$ is a face of both, 
$\lambda$ and~$\tau$. 
Note that we allow here 
a fan to have cones containing
lines.

\begin{theorem}
\label{maina}
Let $\KK$ be an algebraically closed field, 
$M$ a lattice, and $A$ a finitely generated, integral,
$M$-graded $\KK$-algebra with GIT-fan $\Lambda(A)$.
\begin{enumerate}
\item
If $u,v \in \omega(A) \cap M$ is a generating pair,
then the weights $u,v$ lie in a common  GIT-cone 
$\lambda \in \Lambda(A)$.
\item
If $u,v \in \omega(A) \cap M$ lie in a common
GIT-cone $\lambda \in \Lambda(A)$ and $u$ belongs 
to the relative interior 
$\lambda^\circ \subseteq \lambda$,
then $u,v$ is a generating pair.
\end{enumerate}
\end{theorem}

If two weights $u,v \in \omega(A) \cap M$
lie on the boundary of a common GIT-cone
$\lambda \in \Lambda(A)$, 
then no general statement in terms 
of the GIT-fan is possible:
it may happen that $u,v$ 
is generating, and also 
it may happen that $u,v$ is not
generating.
For the first case there are obvious 
examples, and for the latter we 
present the following one.

\begin{example}
Consider the polynomial 
ring 
$A := \KK[T_1,T_2,T_3,T_4]$
over any field~$\KK$.
Then one may define a
$\ZZ^2$-grading of $A$ 
by setting
$$ 
\deg(T_1) \ := \ (4,1),
\quad
\deg(T_2) \ := \ (2,1),
\quad
\deg(T_3) \ := \ (1,2),
\quad
\deg(T_4) \ := \ (1,3).
$$
Any cone in $\QQ^2$ generated by 
a collection of these weights is 
actually an orbit cone, and 
the associated GIT-fan 
looks as follows.
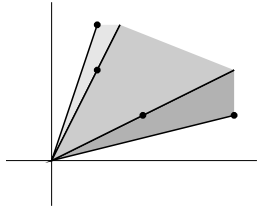
\begin{center}
\begin{picture}(0,0)%
\includegraphics{gitfan.pstex}%
\end{picture}%
\setlength{\unitlength}{1243sp}%
\begingroup\makeatletter\ifx\SetFigFont\undefined%
\gdef\SetFigFont#1#2#3#4#5{%
  \reset@font\fontsize{#1}{#2pt}%
  \fontfamily{#3}\fontseries{#4}\fontshape{#5}%
  \selectfont}%
\fi\endgroup%
\begin{picture}(4974,4074)(1339,-3673)
\end{picture}%

\end{center}
The pair $u := (2,1)$ and $v := (1,2)$
is contained 
in a common GIT-cone but it is not generating:
one directly checks that the 
monomials
$T_1T_2^{n-2}T_3^{n-1}T_4 \in A_{n(u+v)}$ 
can never be
obtained by multiplying elements from
$A_{nu}$ and $A_{nv}$. 
\end{example}

\begin{remark}
In order to compute the 
GIT-fan for concrete examples,
one needs to know the orbit 
cones. Here comes a general 
recipe.

Let $A$ be given by 
homogeneous generators 
and relations,
i.e., we have a graded 
epimorphism
$\KK[T_1, \ldots, T_r] \to A$
and generators 
$q_1, \ldots, q_s$
for its kernel.
With $w_i := \deg(T_i)$,
the orbit cones are
$\cone(w_i; \; i \in I)$,
where $I \subseteq \{1, \ldots, r\}$
satisfies
$$ 
\prod_{i \in I} T_i
\ \not\in \ 
\sqrt{\bangle{q_1^I, \ldots, q_s^I}},
\qquad
\text{with}
\quad
q_j^I 
\ :=\ 
q_j(S_1, \ldots, S_r),
\quad 
S_l
\ := \ 
\begin{cases}
T_l & l \in I,
\\
0 & l \not\in I.
\end{cases}
$$
So, finding the sets of weights 
generating an orbit cone, amounts
to testing for radical ideal 
membership,
which can be performed quite 
efficiently by appropriate
computer algebra systems.
\end{remark}

\begin{remark}
For the polynomial ring
$A = \KK[T_1,\ldots,T_r]$,
the property of being a generating
pair can be formulated as follows 
in a purely combinatorial manner.

Let the grading arise from a linear 
map $Q \colon \ZZ^r \to M$, 
$e_i \mapsto \deg(T_i)$.
Then the weight cone $\omega(A)$ 
is the $Q$-image of the positive 
orthant $\gamma \subseteq \QQ^r$,
and for any integral 
$u \in \omega(A)$, we have the 
polyhedron 
$\Delta_u := Q^{-1}(u) \cap \gamma$.
A pair $u,v \in \omega(A) \cap M$ 
is generating if and only if
there exists an $m > 0$ such that for any 
$k > 0$ one has
\begin{eqnarray*}
\left( \Delta_{kmu} \cap \ZZ^r \right)
\ + \
\left( \Delta_{kmv} \cap \ZZ^r \right)
& = & 
\Delta_{km(u+v)} \ \cap \ \ZZ^r.
\end{eqnarray*}
\end{remark}

In order to present the second result, 
we have to recall from~\cite[Sec.~2]{BeHa2} some more facts 
concerning the GIT-fan.
For any $u \in \omega(A) \cap M$, we have 
an associated nonempty set of semistable 
points: 
$$
X(u)
\ := \ 
\bigcup_{\substack{f \in A_{nu}, \\ n > 0}} X_f
\ = \ 
\{x \in X; \; u \in \omega(x) \}.
$$
We have $X(u) \subseteq X(v)$ 
if and only if the GIT-cone $\lambda(v)$ 
is a face of $\lambda(u)$.
In particular, $u,v \in \omega(A) \cap M$ 
define the same set of semistable points
if and only if they belong to the relative 
interior of a common GIT-cone.

Each set of semistable 
points $X(u)$
admits a good quotient
$X(u) \to Y(u)$
for the action of $T$.
For
$X(u) \subseteq X(v)$,
there is an induced projective 
morphism 
$Y(u) \to Y(v)$
of the quotient spaces.
In particular, if $u,v$ lie in a
common GIT-cone, then we obtain 
a commutative diagram
\begin{equation}
\label{uvdiag}
\vcenter{
\xymatrix{
& 
Y(u+v)
\ar[dl]_{\kappa_u}
\ar[dr]^{\kappa_v}
\ar[dd]_{\kappa}
&
\\
Y(u)
&
&
Y(v)
\\
&
Y(u) \times Y(v)
\ar[ul]^{\pi_u}
\ar[ur]_{\pi_v}
&
}
}
\end{equation}
We denote the image of 
the downwards map $\kappa$ 
by $Z(u,v) := \kappa(Y(u + v))$.
Moreover, we consider the (open) set
$W(A) := \{x \in X; \; \omega(x) = \omega(A)\}$
of points having a generic orbit cone.
For a factorial $A$, 
we then obtain the following
characterization of the 
generating property for 
a pair $u,v$ in the 
relative interior
$\omega(A)^\circ$
of $\omega(A)$.

\begin{theorem}
\label{mainb}
Let $\KK$, $M$ and $A$ be as in~\ref{maina}.
Moreover, suppose that
$A$ is factorial and that
$X \setminus W(A)$ is of codimension
at least two in $X$.
Then, for any two $u,v \in \omega(A)^\circ$
belonging to a common GIT-cone,
the following statements are 
equivalent.
\begin{enumerate}
\item
The pair $u,v$ is generating.
\item
The variety $Z(u,v)$ is normal.
\end{enumerate}
\end{theorem}

\begin{remark}
Under slightly sharper conditions on
the algebra $A$ as posed in Theorem~\ref{mainb}, 
one may view $A$ as the ``Cox ring''
of certain varieties, see~\cite{BeHa1}.
Theorem~\ref{mainb} then tells about 
surjectivity properties of the multiplication
map for global sections of divisors.
\end{remark}

We would like to thank D.~Timashev for fruitful 
discussions 
and the referee for careful reading and
helpful comments.  

\section{Proof of the results}

The setup is the same as in the first section.
In particular, $M$ is a lattice, and
$A$ is a finitely generated, integral 
algebra over an 
algebraically closed field $\KK$.
We consider again the corresponding
affine variety $X := \Spec(A)$,
and the action of the torus $T := \Spec(\KK[M])$
on $X$ defined by the $M$-grading of $A$.

In a first step, we give a more algebraic 
characterization of the GIT-fan.
For $u,v \in \omega(A) \cap M$,
we will work in terms of the following
subalgebras: 
$$ 
A(u)
\ := \
\bigoplus_{n \in \ZZ_{\ge 0}} A_{nu},
\qquad
A(u,v)
\ := \ 
\bigoplus_{n \in \ZZ_{\ge 0}} A_{nu} \mal A_{nv}.
$$
Clearly, $A(u,v)$ is contained in $A(u+v)$. 
We call $A(u,v)$ {\em large\/} in $A(u+v)$,
if the ideals $A(u,v)_+ \subseteq A(u,v)$ 
and $A(u+v)_+ \subseteq A(u+v)$ generated
by the homogeneous parts of strictly positive degree
satisfy
$$
\sqrt{ \bangle{A(u,v)_+} }
\ = \ 
A(u+v)_+
\ \subseteq \ 
A(u+v).
$$

\begin{proposition}
\label{largesubalg}
Let $M$ be a lattice, and $A$ 
an $M$-graded, finitely generated,
integral $\KK$-algebra.
Then, for any two $u,v \in \omega(A)$,
the following statements are equivalent.
\begin{enumerate}
\item
There is a GIT-cone $\lambda \in \Lambda$
satisfying $u,v \in \lambda$.
\item
We have $X(u) \cap X(v) = X(u+v)$.
\item
The algebra $A(u,v)$ is large in $A(u+v)$.
\end{enumerate}
\end{proposition}

\begin{proof}
We begin with the equivalence of (i) and~(ii). 
If~(i) holds, then every orbit cone $\omega(x)$ 
containing $u+v$ must contain $u$ and $v$ as well.
This gives
\begin{eqnarray*}
x \ \in \ X(u) \cap X(v)
& \iff & 
u, v \ \in \ \omega(x)
\\ 
& \iff & 
u + v \ \in \ \omega(x)
\\ 
& \iff & 
x \in \ X(u+v).
\end{eqnarray*}
Conversely, if (ii)~holds, then
we see that $\lambda(u)$ and $\lambda(v)$ 
are faces of $\lambda(u+v)$. 
Thus, we have $u,v \in \lambda(u+v)$.

For the equivalence of (ii) and (iii)
note that for any $w \in \omega(A) \cap M$
the complement $X \setminus X(w)$
equals the zero set $V(A(w)_+)$.
Thus, setting $w := u+v$, we obtain
\begin{eqnarray*}
X(u) \cap X(v)
\ = \ 
X(w)
& \iff & 
V(A(u)_+) \cup V(A(v)_+)
\ = \ 
V(A(w)_+)
\\
& \iff & 
V(A(u)_+ \cdot A(v)_+)
\ = \ 
V(A(w)_+).
\end{eqnarray*}
The latter property holds if and only if
the ideals generated by 
$A(u)_+ \cdot A(v)_+$ and $A(w)_+$
have the same radical in $A$.
This holds if and only if
they generate the same radical ideal
in $A(w)$,
which eventually is equivalent to 
$A(u,v)$ being a large 
subalgebra of $A(w)$.
\end{proof}

This observation enables us to decide
whether or not two weights $u,v$ 
belong to a common GIT-cone by just 
looking at $A(u)$, $A(v)$ and $A(u+v)$.
As a consequence, we may produce examples
of nontrivial affine varieties with
simple variation of GIT-quotients.

Recall that a point $x \in X(u)$ in a set 
$X(u) \subseteq X$ of semistable points
is said to be {\em stable\/}, 
if its orbit $T \mal x$ is closed in $X(u)$
and of maximal dimension.
If the set $X(u)$ consists of stable points,
then the fibres of the quotient map 
$X(u) \to Y(u)$ are precisely the 
$T$-orbits of $X(u)$.

\begin{corollary}
Let $M$ be a lattice, and let $A$ 
be an $M$-graded, finitely generated,
integral $\KK$-algebra.
Given $\lambda \in \Lambda(A)$,
consider the (finitely generated) 
algebra
\begin{eqnarray*}
A'
& := & 
\bigoplus_{u \in \lambda \cap M}
A_u.
\end{eqnarray*}
Then the corresponding action of 
the torus $T = \Spec(\KK[M])$ on 
the affine variety $X' = \Spec(A')$ 
has the following properties.
\begin{enumerate}
\item
The GIT-fan $\Lambda(A')$ associated 
to $A'$
is the fan of faces of the cone 
$\lambda \in \Lambda(A)$.
\item
The union $W \subseteq X'$ of all 
$T$-orbits of maximal dimension
is a set of semistable points, 
and every $x \in W$ is stable.
\end{enumerate}
\end{corollary}

\begin{proof}
To see~(i), note first that $\Lambda(A')$ 
subdivides $\omega(A')=\lambda$. Moreover, 
Proposition~\ref{largesubalg}~(iii) 
implies that two
weights $u,v \in \lambda$ lie in 
a common cone of $\Lambda(A')$ 
if and only if they lie in 
a common cone of $\Lambda(A)$. 

For~(ii), note that the dimension of 
an orbit cone $\omega(x)$ equals that 
of the orbit $\dim(T \cdot x)$. 
Since $\lambda \in \Lambda(A')$ is the
only cone of maximal dimension, we 
obtain 
$$
W \ = \ 
\{x \in X; \; \omega(x) = \lambda\}
\ = \ X'(u)
$$ 
for any $u$ from the relative interior 
of $\lambda$.
Since all orbits in $W$ have the same 
dimension, each of them is closed in $W$.
\end{proof}

The next step is a geometric characterization
of the GIT-fan.
It is given in terms of the map
$\kappa \colon Y(u+v) \to Y(u) \times Y(v)$ 
introduced in the diagram~\ref{uvdiag}.

\begin{proposition}
\label{gengeomchar1}
Let $u,v \in \omega(A) \cap M$
belong to a common GIT-cone
$\lambda \in \Lambda(A)$.
Then, in the setting of~\ref{uvdiag}, 
the following statements are
equivalent:
\begin{enumerate}
\item
The pair $u,v \in \omega(A) \cap M$ 
is generating.
\item
The map $\kappa \colon Y(u+v) \to Y(u) \times Y(v)$ 
is a closed embedding.
\end{enumerate}  
\end{proposition}

\begin{proof}
Recall that the quotient spaces
$Y(w) = \Proj(A(w))$ are projective over 
$Y_0 = \Spec(A_0)$.
Moreover, denoting by 
$q \colon X(w) \to Y(w)$ the quotient
map, we obtain for $n \in \ZZ_{\ge 0}$
a sheaf on $Y(w)$, namely
$$
\mathcal{L}_{nw}
\ := \
\left( q_* \mathcal{O}_{X(w)} \right)_{nw}
\ = \ 
\mathcal{O}_{Y(w)}(n).
$$

Replacing $u$ with a large 
multiple,
we may assume that $A(u)$ 
is generated as an $A_0$-algebra
by the component $A_{u}$, 
and that for any $n \in \ZZ_{\ge 1}$
the canonical maps
\begin{eqnarray*}
\imath_{nu} \colon
A_{nu} 
& \to &
\Gamma(Y(u),\mathcal{L}_{nu})
\end{eqnarray*}
are surjective,
see~\cite[Exercise~II.5.9]{Har}.
Note that then
$\mathcal{L}_u$ is 
an ample invertible  
sheaf on~$Y(u)$.
Of course, we may arrange the 
same situation for $v$ and 
$u+v$.

On $Y(u) \times Y(v)$ we have the ample 
invertible sheaves
$\mathcal{E}_n := \pi_u^* \mathcal{L}_{nu} \otimes \pi_v^* \mathcal{L}_{nv}$.
We claim that the natural map
\begin{eqnarray*}
\Gamma(Y(u),\mathcal{L}_{nu}) \otimes \Gamma(Y(v), \mathcal{L}_{nv})
& \to & 
\Gamma(Y(u) \times Y(v), \mathcal{E}_n)
\end{eqnarray*}
is an isomorphism. Indeed, using the projection 
formula, we obtain canonical isomorphisms
$$
\Gamma(Y(u)\times Y(v), \mathcal{E}_n)
\ \cong \ 
\Gamma(Y(u),{\pi_u}_* \mathcal{E}_n)
\ \cong \
\Gamma(Y(u), \mathcal{L}_{nu} \otimes {\pi_u}_* \pi_v^*\mathcal{L}_{nv}).
$$
We look a bit closer at ${\pi_u}_* \pi_v^*\mathcal{L}_{nv}$.
Given an open subset $U \subseteq Y(u)$, we denote 
 by $\pi^U_v \colon U \times Y(v) \to Y(v)$
the restricted projection. Then we have
$$
\Gamma(U,{\pi_u}_*\pi_v^*\mathcal{L}_{nv})
\ = \
\Gamma(U\times Y(v),\pi_v^*\mathcal{L}_{nv})
\ \cong \ 
\Gamma(Y(v), \mathcal{L}_{nv} 
\otimes
{\pi^U_v}_* \mathcal{O}_{U\times Y(v)}).
$$
Likewise, one obtains 
${\pi^U_v}_* \mathcal{O}_{U\times Y(v)} \cong
\Gamma(U, \mathcal{O}_U) \otimes \mathcal{O}_{Y(v)}$
for any affine open set $U \subseteq Y(u)$.
Consequently, we have a canonical isomorphism
$$
\Gamma(U,{\pi_u}_* \pi_v^*\mathcal{L}_{nv})
\ \cong \
\Gamma(U,\mathcal{O}_U) \otimes \Gamma(Y(v),\mathcal{L}_{nv}).
$$
This in turn shows
${\pi_u}_* \pi_v^*\mathcal{L}_{nv}
\cong
\mathcal{O}_{Y(u)}\otimes\Gamma(Y(v),\mathcal{L}_{nv})$,
and our claim follows.
Thus, we arrive at a commutative diagram
$$ 
\xymatrix{
A_{nu} \otimes A_{nv}
\ar[rr]^{\mu_{n}}
\ar[d]_{\cong}
& &
A_{nu+nv}
\ar[d]^{\cong}
\\
{\Gamma(Y(u) \times Y(v), \mathcal{E}_n)}
\ar[rr]_{\kappa_n^*}
& &
{\Gamma(Y(u+v), \mathcal{L}_{nu+nv})}
}
$$
where the upper horizontal arrow is the 
multiplication map we are interested in,
and the lower horizontal arrow is the 
canonical pullback map
\begin{eqnarray*}
\kappa_{n}^* \colon
\Gamma(Y(u) \times Y(v),
\mathcal{E}_n)
& \to &
\Gamma(Y(u+v), \mathcal{L}_{nu+nv})
\\
\pi_u^* f \otimes \pi_v^* g 
& \mapsto &
\kappa_u^* f  \cdot \kappa_v^* g.
\end{eqnarray*}

Now, note that the morphism 
$\kappa \colon Y(u+v) \to Y(u) \times Y(u)$ 
is induced from the multiplication map,
because we have 
$$ 
Y(u) \times Y(v)
\ = \ 
\Proj\left( \bigoplus_{n \ge 0} A_{nu} \otimes A_{nv} \right),
\qquad
Y(u+v)
\ = \ 
\Proj \left( \bigoplus_{n \ge 0} A_{nu + nv} \right).
$$
Thus, the assertion follows from the basic fact that
$\kappa$ is a closed embedding if and only if
there is an $l > 1$ such that
$\mu_{ln}$ are surjective for any $n>0$.
\end{proof}

\begin{proof}[Proof of Theorem~\ref{maina}]
If $u,v \in \omega \cap M$ is a generating pair,
then the algebra $A(u,v)$ is large in~$A(u+v)$.
Thus, the first assertion follows from 
Proposition~\ref{largesubalg}.
To see the second one, note that both, $u$ and 
$u +v$, lie in  the relative interior
$\lambda^{\circ}$ of the GIT-cone 
$\lambda \in \Lambda(A)$.
Thus, $Y(u+v) \to Y(u)$ is an isomorphism,
and the statement follows from 
Proposition~\ref{gengeomchar1}.
\end{proof}

\begin{proof}[Proof of Theorem~\ref{mainb}]
First note that the set $W := W(A) \subseteq X$ 
consisting of all $x \in X$ with 
orbit cone $\omega(x) = \omega(A)$ admits
a geometric quotient $V := W/T$ and that 
for any $w \in \omega(A)^\circ$,
the inclusion $W \subseteq X(w)$ 
induces an open embedding 
$V \to Y(w)$
of the quotient spaces.
Since $W \subseteq X$ has 
a complement of codimension 
at least two in $X$,
the same must hold for
the image of $V$ in $Y(w)$.
Moreover, as a good quotient
space of a normal variety,
$Y(w)$ is normal.
Thus, 
$V \to Y(w)$ 
is a $V$-embedding 
in the sense of~\cite[Sec.~2]{ArHa}.

To proceed, consider the morphisms
of~\ref{uvdiag}.
Clearly,
$\kappa_u \colon Y(u+v) \to Y(u)$ 
and $\kappa_v \colon Y(u+v) \to Y(v)$
are morphisms of $V$-embeddings, 
that means that
we have a commutative diagram
$$ 
\xymatrix{
& 
V 
\ar[d]
\ar[dr]
\ar[dl]
\\
Y(u)
&
Y(u+v)
\ar[l]^{\kappa_u}
\ar[r]_{\kappa_v}
&
Y(v)
}
$$
Now consider the map
$\kappa \colon Y(u+v) \to Y(u) \times Y(v)$ 
of~\ref{uvdiag}, and denote 
its image by $Z := Z(u,v)$.
Then $\kappa$ lifts to the normalization
$Z' \to Z$, and we obtain a 
commutative diagram
$$
\xymatrix{
& 
Y(u+v)
\ar[dl]_{\kappa_u}
\ar[dr]^{\kappa_v}
\ar[d]
&
\\
Y(u)
&
Z'
\ar[l]
\ar[r]
\ar[d]
&
Y(v)
\\
&
Z
\ar[ul]^{\pi_u}
\ar[ur]_{\pi_v}
&
}
$$

Lifting 
$V \to Y(u+v) \to Z$ to $Z'$ defines a 
$V$-embedding $V \to Z'$.
According to~\cite[Prop.~2.3]{ArHa}, 
there is an open $T$-invariant
subset $W' \subseteq X$ with
good quotient $W' \to Z'$ by the
$T$-action such that
$V \to Z'$ is induced by 
the inclusion $W \subseteq W'$.

Moreover, the map
$Y(u+v) \to Z'$ as well as the maps
$Z' \to Y(u)$ and $Z' \to Y(v)$ 
are morphisms of $V$-embeddings.
Thus, \cite[Prop.~2.4]{ArHa} tells us
that they are induced by inclusions 
of sets of semistable points
$$
X(u+v) \ \subseteq \ W',
\qquad
W' \ \subseteq \ X(u),
\qquad
W' \ \subseteq \ X(v).
$$
By Proposition~\ref{largesubalg}, we have 
$X(u+v) = X(u) \cap X(v)$.
This shows $W' = X(u+v)$. Thus, 
the map $Y(u+v) \to Z'$ is an isomorphism.
{From} this we see that the map
$\kappa \colon Y(u+v) \to Y(u) \times Y(v)$ 
is a closed 
embedding if and only if $Z$ is normal.
The assertion then follows from
Proposition~\ref{gengeomchar1}.
\end{proof}

\end{document}